\documentstyle[12pt]{amsart}
\newcommand{\g}{\goth}
\newcommand{\gtg}{\mbox{\g g}}

\newcommand{\gtsl}{\mbox{\g sl}}
\newcommand{\gtgl}{\mbox{\g gl}}
\newcommand{\hgtsl}{\mbox{$\hat{\gtsl}$}}
\newcommand{\hgtgl}{\mbox{$\hat{\gtgl}$}}

\newcommand{\gth}{\mbox{\g h}}

\newcommand{\nc}{\mbox{${\Bbb C}$}}
\newcommand{\nz}{\mbox{${\Bbb Z}$}}

\newcommand{\cA}{\mbox{${\cal A}$}}

\newcommand{\cD}{\mbox{${\cal D}$}}
\newcommand{\cF}{\mbox{${\cal F}$}}

\newcommand{\cP}{\mbox{${\cal P}$}}
\newcommand{\cR}{\mbox{${\cal R}$}}
\newcommand{\id}{\mbox{\rm id}}

\newcommand{\End}{\mbox{\rm End}}

\newcommand{\vep}{\varepsilon}

\theoremstyle{plain}
 \newtheorem{thm}{Theorem}[section]
 
 \newtheorem{lemma}[thm]{Lemma}
 \newtheorem{cor}[thm]{Corollary}
\theoremstyle{definition}
 \newtheorem{defn}{Definition}[section]

\theoremstyle{remark}

 \newtheorem{ack}{Acknowledgement}

\begin{document}
\title{A central extension of $\cD Y_{\hbar}(\gtgl_2)$ \\ 
      and its vertex representations}
\author{Kenji Iohara$^\dagger$}
\address[K. Iohara]{Department of Mathematics, Faculty of Science,
                    Kyoto University, Kyoto 606, Japan.}
\author{Mika Kohno}
\address[M. Kohno]{Department of Mathematical Sciences, 
                 University of Tokyo, Tokyo 153, Japan.}
\thanks{$\dagger$:JSPS Research Fellow}
\date{\today}
\maketitle
\begin{abstract}

A central extension of $\cD Y_{\hbar}(\gtgl_2)$ is proposed. 
The bosonization of level $1$ module and vertex operators 
are also given.
\end{abstract}
\pagestyle{plain}
\section{Introduction}\hspace{1 in} \\

 A quantum algebra called Yangian was introduced in \cite{dr_1} in connection
 with the rational solutions of quantum Yang-Baxter equation.
The Yangian double $\cD Y_{\hbar}(\gtg)$ is a quantum double \cite{dr_icm}
of the Yangian. This algebra has been introduced in the literature
in terms of Chevalley generators \cite{leclair}, $T^{\pm}$-matrix
 \cite{bern} and  Drinfel'd generators \cite{tolstoy}. It has also 
applications to the massive
integrable field theory \cite{leclair},\cite{smirnov},\cite{bern}.
In these works, $\cD Y_{\hbar}(\gtg)$ is defined as a deformation of the 
enveloping algebra of the loop algebra $\gtg[t,t^{-1}]$, without central
extension. In view of the development in lattice models \cite{jimbo_miwa}, 
it seems natural and useful to construct an algebra having highest 
weight modules and vertex operators. The purpose of this article is to
define a central extension of Yangian double $\cD Y_{\hbar}(\gtg)$
for $\gtg=\gtgl_2,\gtsl_2$, which we again denote by the same symbol.
We also present analogs of level 1 modules and vertex operators 
in bosonized form.

 We introduce $\cD Y_{\hbar}(\gtgl_2)$ by means of
 the quantum inverse scattering method, or `$RTT$ formalism' 
\cite{bern},\cite{resh},\cite{QISM}. There is another important
set of generators called Drinfel'd generators. The correspondence 
between $T^{\pm}$-matrix and the Drinfel'd generators are known for
$Y_{\hbar}(\gtsl_n)$ by \cite{dr_2} and for quantized affine algebra
$U_q(\hgtgl_n)$ by \cite{ding}. Following the original argument 
\cite{dr_2}, we obtain the Drinfel'd generators 
of $\cD Y_{\hbar}(\gtg)$ for $\gtg=\gtgl_2,\gtsl_2$ in terms of
currents (Theorem \ref{gl_2},Corollary \ref{sl_2}). This recovers 
the result obtained in \cite{tolstoy} at level $0$. 
 We expect that such structures as highest weight modules and vertex
operators for quantized affine algebras have natural counterpart for
our algebra $\cD Y_{\hbar}(\gtg)$. 
Unfortunately, lacking an appropriate definition of the 
grading operator $d$, we still do not know how to define highest 
weight modules. Nevertheless, with the aid of the Drinfel'd generators, 
we can construct explicitly an analogue 
of level $1$ module $\cF_{i,s}$ in terms of bosons  
(theorem \ref{level 1}) in a way parallel to the
 quantum affine version \cite{jing}. 
We also give an analogue of vertex operators
(Theorem \ref{VOp}) acting on $\cF_{i,s}$, as certain intertwiners
of $\cD Y_{\hbar}(\gtg)$-module and work out their
commutation relations (Theorem \ref{comm}). 

 The text is organized as follows. In section 2 we define 
$\cD Y_{\hbar}(\gtg)$ for $\gtg=\gtgl_2,\gtsl_2$ and gives
its Drinfel'd generators.
In section 3 we construct bosonization of level 1 modules and vertex 
operators. Section 4 contains discussions and remarks.
A detailed account including
the higher rank case will appear elsewhere \cite{kohno}. 

\section{The algebra $\cD Y_{\hbar}(\gtgl_2)$}\hspace{1 in} \\
At level $0$, there are several approaches to define the Yangian double 
$\cD Y_{\hbar}(\gtsl_2)$ \cite{leclair},\cite{bern},\cite{tolstoy}. 
Here we define a central extension of $\cD Y_{\hbar}(\gtg)$ for 
$\gtg=\gtgl_2,\gtsl_2$ following the method of \cite{resh}.

Let $V$ be a $2$-dimensional vector space and  $\cP \in \End(V\otimes V)$
be a permutation operator $\cP v\otimes w =w\otimes v~(v,w\in V)$.
Consider Yang's $R$-matrix normalized as
\[ R(u)=\frac{1}{u+\hbar}\left( uI+\hbar\cP \right) \in
        \End(V\otimes V),\] 
where $\hbar$ is a formal variable.
This $R$-matrix satisfies the following properties:
\begin{description}
\item[Yang-Baxter equation]
\[ R_{12}(u-v)R_{13}(u)R_{23}(v)=
   R_{23}(v)R_{13}(u)R_{12}(u-v), \]
\item[Unitarity]
\[ R_{12}(u)R_{21}(-u)=\id . \]
\end{description}

Here, if $R(u)=\sum a_i \otimes b_i$ with $a_i,b_i \in \End (V)$, then 
$R_{21}(u)=\sum b_i \otimes a_i,~R_{13}(u)=\sum a_i \otimes 1 \otimes b_i$ 
etc.. Set $\cA=\nc[[\hbar]]$.

\begin{defn}
$\cD Y_{\hbar}(\gtgl_2)$ is a topological Hopf algebra over
$\cA$ generated by
$\{ t_{ij}^k |1\leq i,j \leq 2~,~k\in \nz \}$ and $c$. In terms of matrix
generating series
\begin{align*}
&  T^{\pm}(u)=(t_{ij}^{\pm}(u))_{1\leq i,j \leq 2}, \\
&  t_{ij}^{+}(u)=\delta_{ij}
                -\hbar \sum_{k\in \nz_{\geq 0}}t_{ij}^{k}u^{-k-1}~,~
   t_{ij}^{-}(u)=\delta_{ij}
                +\hbar \sum_{k\in \nz_{<0}}t_{ij}^{k}u^{-k-1}, 
\end{align*}
the defining relations are given as follows:
\begin{align*}
[T^{\pm}(u),c]& =0, \\
R(u-v)\overset{1}{T^{\pm}}(u)\overset{2}{T^{\pm}}(v)&=
\overset{2}{T^{\pm}}(v)\overset{1}{T^{\pm}}(u)R(u-v), \\
R(u_{-}-v_{+})\overset{1}{T^{+}}(u)\overset{2}{T^{-}}(v)&=
\overset{2}{T^{-}}(v)\overset{1}{T^{+}}(u)R(u_{+}-v_{-}).
\end{align*}
Here
\[   \overset{1}{T}(u)=T(u)\otimes \id~,
    ~\overset{2}{T}(u)=\id \otimes T(u), \]
$u_{\pm}=u\pm \frac{1}{4}\hbar c$ and similarly for $v$.
Its coalgebra structure is defined as 
\begin{align*}
&  \Delta(t_{ij}^{\pm}(u))=\sum_{k=1}^2
   t_{kj}^{\pm}(u\pm \frac{1}{4}\hbar c_2)\otimes
   t_{ik}^{\pm}(u\mp \frac{1}{4}\hbar c_1), \\             
&  \vep(T^{\pm}(u))=I~,~S({}^tT^{\pm}(u))=[{}^tT^{\pm}(u)]^{-1}, \\
&  \Delta(c)=c\otimes 1+1 \otimes c~,~\vep(c)=0~,~S(c)=-c,
\end{align*}
where $c_1=c\otimes 1$ and $c_2=1\otimes c$.
\end{defn}
Note that the subalgebra generated by
$\{ t_{ij}^k|1\leq i,j\leq 2~,~k\in \nz_{\geq 0} \}$ is  
$Y_{\hbar}(\gtgl_2)$~ \cite{dr_1,dr_2}.  
 
 We introduce the Drinfel'd generators
of $\cD Y_{\hbar}(\gtgl_2)$ exactly in the same way as in
\cite{dr_2}. Namely we define
\begin{align*}
&  X^{+}(u)=t_{11}^{+}(u_{-})^{-1}t_{12}^{+}(u_{-})-
            t_{11}^{-}(u_{+})^{-1}t_{12}^{-}(u_{+}), \\
&  X^{-}(u)=t_{21}^{+}(u_{+})t_{11}^{+}(u_{+})^{-1}-
            t_{21}^{-}(u_{-})t_{11}^{-}(u_{-})^{-1}, \\
&  k_{1}^{\pm}(u)=t_{11}^{\pm}(u),~
   k_{2}^{\pm}(u)=t_{22}^{\pm}(u)-
                  t_{21}^{\pm}(u)t_{11}^{\pm}(u)^{-1}t_{12}^{\pm}(u).
\end{align*}
Notice that $t_{ii}^{\pm}(u)$ have unique inverses as power series
in $\hbar$. The commutation relations among 
$X^{\pm}(u),k_{1}^{\pm}(u),k_{2}^{\pm}(u)$ can be written as follows.
\begin{thm} \label{gl_2}
\begin{align*}
              & k_{i}^{\pm}(u)k_{j}^{\pm}(v)
              =k_{j}^{\pm}(v)k_{i}^{\pm}(u),~
               k_{i}^{+}(u)k_{i}^{-}(v)=k_{i}^{-}(v)k_{i}^{+}(u), \\
              & \frac{u_{\mp}-v_{\pm}}{u_{\mp}-v_{\pm}+\hbar}
               k_{2}^{\mp}(v)^{-1}k_{1}^{\pm}(u)=
               \frac{u_{\pm}-v_{\mp}}{u_{\pm}-v_{\mp}+\hbar}
               k_{1}^{\pm}(u)k_{2}^{\mp}(v)^{-1}, \\
             & \begin{cases}
               k_{1}^{\pm}(u)^{-1}X^{+}(v)k_{1}^{\pm}(u)=
               \dfrac{u_{\pm}-v+\hbar}{u_{\pm}-v}X^{+}(v), &  \\
               k_{1}^{\pm}(u)X^{-}(v)k_{1}^{\pm}(u)^{-1}=
               \dfrac{u_{\mp}-v+\hbar}{u_{\mp}-v}X^{-}(v), &
               \end{cases} \\
             & \begin{cases}
               k_{2}^{\pm}(u)^{-1}X^{+}(v)k_{2}^{\pm}(u)=
               \dfrac{u_{\pm}-v-\hbar}{u_{\pm}-v}X^{+}(v), & \\
               k_{2}^{\pm}(u)X^{-}(v)k_{2}^{\pm}(u)^{-1}=
               \dfrac{u_{\mp}-v-\hbar}{u_{\mp}-v}X^{-}(v), &
               \end{cases} \\
             & (u-v\mp \hbar)X^{\pm}(u)X^{\pm}(v)=
               (u-v\pm \hbar)X^{\pm}(v)X^{\pm}(u), \\
             & [X^{+}(u),X^{-}(v)]=\hbar 
               \left\{ \delta(u_{-}-v_{+})
                       k_{2}^+(u_{-})k_{1}^+(u_{-})^{-1}-
                       \delta(u_{+}-v_{-})
                       k_{2}^-(v_{-})k_{1}^-(v_{-})^{-1} \right\}.
\end{align*}
here $\delta(u-v)=\sum_{k\in \nz} u^{-k-1}v^k$.
\end{thm}
 To decompose $\cD Y_{\hbar}(\gtgl_2)$
into two subalgebras 
$\cD Y_{\hbar}(\gtsl_2)$ and a Heisenberg subalgebra,
we introduce the following currents:
\begin{align*}
&  H^{\pm}(u)=k_{2}^{\pm}(u+\frac{1}{2}\hbar )
              k_{1}^{\pm}(u+\frac{1}{2}\hbar )^{-1}~,~
   K^{\pm}(u)=k_{1}^{\pm}(u-\frac{1}{2}\hbar)
              k_{2}^{\pm}(u+\frac{1}{2}\hbar), \\
&  E(u)=\frac{1}{\hbar}X^{+}(u+\frac{1}{2}\hbar )~,~
   F(u)=\frac{1}{\hbar}X^{-}(u+\frac{1}{2}\hbar ). 
\end{align*}
We define $\cD Y_{\hbar}(\gtsl_2)$ to be the subalgebra of
 $\cD Y_{\hbar}(\gtgl_2)$ generated by
$H^{\pm}(u)~,~E(u)~,~ F(u)~$ and $c$. A Heisenberg subalgebra 
of $ \cD Y_{\hbar}(\gtgl_2) $ generated by $K^{\pm}(u)$ 
commute with all of the elements of $ \cD Y_{\hbar}(\gtsl_2) $.
In terms of these generators, the above commutation relations can be 
rephrased as follows.
\begin{cor}
\begin{align*}\label{sl_2}
& [H^{\pm}(u),H^{\pm}(v)]=0, \\
&  (u_{\mp}-v_{\pm}+\hbar)(u_{\pm}-v_{\mp}-\hbar)
  H^{\pm}(u)H^{\mp}(v) \\
& =(u_{\mp}-v_{\pm}-\hbar )(u_{\pm}-v_{\mp}+\hbar )
  H^{\mp}(v)H^{\pm}(u), \\
& [K^{\pm}(u),K^{\pm}(v)]=0, \\
&  \frac{u_{-}-v_{+}-\hbar}{u_{-}-v_{+}+\hbar}K^{+}(u)K^{-}(v)=   
   K^{-}(v)K^{+}(u)\frac{u_{+}-v_{-}-\hbar}{u_{+}-v_{-}+\hbar}, \\
& \begin{cases}
  H^{\pm}(u)^{-1}E(v)H^{\pm}(u)=
  \dfrac{u_{\pm}-v-\hbar}{u_{\pm}-v+\hbar}E(v), & \\
  H^{\pm}(u)F(v)H^{\pm}(u)^{-1}=
  \dfrac{u_{\mp}-v-\hbar}{u_{\mp}-v+\hbar}F(v), &
  \end{cases} \\
& [K^{\sigma}(u),H^{\pm}(v)]=[K^{\sigma}(u),E(v)]
 =[K^{\sigma}(u),F(v)]=0 \quad \forall \sigma=\pm, \\
& (u-v-\hbar )E(u)E(v)=(u-v+\hbar )E(v)E(u), \\
& (u-v+\hbar )F(u)F(v)=(u-v-\hbar )F(v)F(u), \\
& [E(u),F(v)]=\frac{1}{\hbar}
  \left\{ \delta(u_{-}-v_{+})H^{+}(u_{-})-
          \delta(u_{+}-v_{-})H^{-}(v_{-})\right\}.
\end{align*}
\end{cor}

To compare with the known results at $c=0$ \cite{tolstoy},
let us write down the commutation relations componentwise.
 The Fourier components of the generating series
$H^{\pm}(u),E(u),F(u)$ are of the following form:
\begin{align*}
& H^{+}(u)=1+\hbar\sum_{k\geq 0}h_{k}u^{-k-1},~
  H^{-}(u)=1-\hbar\sum_{k<0}h_{k}u^{-k-1}, \\
& E(u)=\sum_{k\in \nz}e_{k}u^{-k-1},~
  F(u)=\sum_{k\in \nz}f_{k}u^{-k-1}.
\end{align*}
For $c=0$, the commutation relations of $ \cD Y_{\hbar}(\gtsl_2)$
in terms of the above Fourier component look simple as follows:
\begin{align*}
& [h_{k},h_{l}]=0,~ [h_{0},x_{l}^{\pm}]=\pm2x_{l}^{\pm},~
  [x_{k}^+,x_{l}^-]=h_{k+l}, \\
& [h_{k+1},x_{l}^{\pm}]-[h_{k},x_{l+1}^{\pm}]=
  \pm \hbar [h_{k},x_{l}^{\pm}]_{+}, \\
& [x_{k+1}^{\pm},x_{l}^{\pm}]-[x_{k}^{\pm},x_{l+1}^{\pm}]=
  \pm \hbar [x_{k}^{\pm},x_{l}^{\pm}]_{+}, 
\end{align*}
for $k,l \in \nz$, where we set $x_{k}^{+}=e_{k}~,~x_{k}^{-}=f_{k}$ and
$[x,y]_{+}=xy+yx$ for $x,y \in \cD Y_{\hbar}(\gtsl_2)$.
These relations are the same as in \cite{tolstoy} with $\hbar=1$.
The set $\{ h_k, x_k^{\pm}| k \in \nz_{\geq 0} \}$ 
provides the Drinfel'd generators of $Y_{\hbar}(\gtsl_2)$ \cite{dr_2}.

Let us set
\[E^{\pm}(u)=\frac{1}{\hbar}t_{11}^{\pm}(u_{\mp}+\frac{1}{2}\hbar)^{-1}
                            t_{12}^{\pm}(u_{\mp}+\frac{1}{2}\hbar),~
  F^{\pm}(u)=\frac{1}{\hbar}t_{21}^{\pm}(u_{\pm}+\frac{1}{2}\hbar)
                            t_{11}^{\pm}(u_{\pm}+\frac{1}{2}\hbar)^{-1}, \]
so that $E(u)=E^{+}(u)-E^{-}(u),~F(u)=F^{+}(u)-F^{-}(u)$. One can calculate
the coproduct of the currents $E^{\pm}(u),F^{\pm}(u),H^{\pm}(u),K^{\pm}(u)$
and the results are the following.
\begin{lemma}\label{coproduct}
\begin{align*}
1)~ \Delta(E^{\pm}(u)) &=
     E^{\pm}(u)\otimes 1 \\ 
   & + \sum_{k \geq 0}
     (-1)^k {\hbar}^{2k}F^{\pm}(u\mp \frac{1}{2}\hbar c +\hbar)^k 
     H^{\pm}(u_{\mp}) \\
   & \quad \otimes E^{\pm}(u\mp \frac{1}{2}\hbar c_1)^{k+1}, \\
2)~ \Delta(F^{\pm}(u)) &=
     1\otimes F^{\pm}(u) \\
   & + \sum_{k \geq 0}
     (-1)^k {\hbar}^{2k}F^{\pm}(u \pm \frac{1}{2}\hbar c_2)^{k+1} \\
   & \quad \otimes
     E^{\pm}(u \pm \frac{1}{2}\hbar c -\hbar)^k
     H^{\pm}(u_{\pm}), \\
3)~ \Delta(H^{\pm}(u)) &=
     \sum_{k \geq 0} (-1)^{k}(k+1){\hbar}^{2k}
     F^{\pm}(u_{\mp}+\hbar \pm \frac{1}{4}\hbar c_2)^k 
     H^{\pm}(u\pm \frac{1}{4}\hbar c_2) \\
   & \quad \otimes 
     E^{\pm}(u_{\pm}-\hbar \mp \frac{1}{4}\hbar c_1)^k
     H^{\pm}(u\mp \frac{1}{4}\hbar c_1), \\
4)~  \Delta(K^{\pm}(u)) &=
     K^{\pm}(u\pm \frac{1}{4}\hbar c_2)\otimes 
     K^{\pm}(u\mp \frac{1}{4}\hbar c_1).
\end{align*}
\end{lemma} 

\section{Bosonization of level $1$ module}\hspace{1 in} \\
We expect that there is an appropriate definition of highest weight
modules which can be regarded as deformation of their affine counterpart.
At the moment, we have difficulty due to the lack of a proper 
definition of the grading operator $d$.
Here we construct level $1$ $ \cD Y_{\hbar}(\gtgl_2)$-module 
and vertex operators directly in terms of bosons. 

Let $\gth=\nc \vep_1 \oplus \nc \vep_2$ be a Cartan subalgebra of
$\gtgl_2$, $\overline{Q}=\nz \alpha~(\alpha=\vep_1-\vep_2)$ 
be the root lattice of $\gtsl_2$,  
$\overline{\Lambda}_i=\Lambda_i-\Lambda_0$ be the classical part of the
$i$-th fundamental weight and $(\cdot, \cdot)$ be the standard bilinear
form defined by $(\vep_i,\vep_j)=\delta_{ij}$.
Let us introduce bosons $\{ a_{i,k}|i=1,2,~k\in \nz 
\setminus \{0 \} \}$ satisfying:
\[ [a_{i,k},a_{j,l}]=k\delta_{i,j}\delta_{k+l,0}. \]
Set
\[ \cF_{i,s}:=\cA [a_{j,-k}(j=1,2,~k\in \nz_{>0})]\otimes
              \cA [\overline{Q}]
              e^{\overline{\Lambda}_i+\frac{s}{2}(\vep_1+\vep_2)}
   \quad (i=0,1), \]
where $s$ is a complex parameter and $\cA[\overline{Q}]$ is 
the group algebra of $\overline{Q}$ over $\cA$. 
On this space, we define the action of the operators 
$a_{j,k},\partial_{\vep_j},e^{\vep_j}$~$(j=1,2)$ by
\begin{align*}
a_{j,k}\cdot f\otimes e^{\beta}& =\begin{cases}
                      a_{j,k}f \otimes e^{\beta}          & k< 0  \\
            \text{$[a_{j,k},f]$}  \otimes e^{\beta} \quad & k> 0
                                  \end{cases}, \\
\partial_{\vep_j}\cdot f\otimes e^{\beta}& 
        =(\vep_j,\beta)f\otimes e^{\beta}, 
         \qquad \text{for}~f\otimes e^{\beta}\in \cF_{i,s} \\
e^{\vep_j}\cdot f\otimes e^{\beta} & 
=f\otimes e^{\vep_j+\beta}.
\end{align*}

\begin{thm} \label{level 1}
The following assignment defines a 
$ \cD Y_{\hbar}(\gtgl_2)$-module structure on $\cF_{i,s}$.
\begin{align*}
k_j^+(u)&\mapsto \exp \left[ -\sum_{k>0}
          \frac{a_{j,k}}{k}\left \{ (u+\frac{1}{2}\hbar)^{-k}
                                  -(u-\frac{1}{2}\hbar)^{-k}
          \right \} \right]
          \left( \frac{u-\frac{1}{2}\hbar}{u+\frac{1}{2}\hbar}
          \right)^{\partial_{\vep_j}}, \\
k_j^-(u)&\mapsto \exp \left[ \sum_{k>0,r<j}\frac{a_{r,-k}}{k}
          \left \{ u^k-(u-\hbar)^k \right \}+
                      \sum_{k>0,r>j}\frac{a_{r,-k}}{k}
          \left \{ (u+\hbar)^k-u^k \right \} \right], \\                  \frac{1}{\hbar}X^{+}(u)& \mapsto \exp \left[ -\sum_{k>0}
          \frac{a_{1,-k}}{k}(u-\frac{3}{4}\hbar)^{k}+
          \sum_{k>0}\frac{a_{2,-k}}{k}(u+\frac{1}{4}\hbar)^{k}
          \right] \\
               & \times \exp \left[ \sum_{k>0}
          \frac{a_{1,k}-a_{2,k}}{k}(u+\frac{1}{4}\hbar)^{-k}
          \right] e^{\alpha}
          \left[ u+\frac{1}{4}\hbar
          \right]^{\partial_{\alpha}}, \\
\frac{1}{\hbar}X^{-}(u)& \mapsto \exp \left[ \sum_{k>0}
          \frac{a_{1,-k}}{k}(u-\frac{1}{4}\hbar)^{k}-
          \sum_{k>0}\frac{a_{2,-k}}{k}(u+\frac{3}{4}\hbar)^{k}
          \right] \\
               & \times \exp \left[ \sum_{k>0}
          \frac{-a_{1,k}+a_{2,k}}{k}(u-\frac{1}{4}\hbar)^{-k}
          \right] e^{-\alpha}
          \left[ u-\frac{1}{4}\hbar
          \right]^{-\partial_{\alpha}}.
\end{align*}
\end{thm}
Similar formulas are established for $ \cD Y_{\hbar}(\gtsl_2)$
by using one boson $a_k~(k\in \nz)$ \cite{kohno}. 
We do not know how to characterize $\cF_{i,s}$ conceptually as
$ \cD Y_{\hbar}(\gtgl_2)$-module.

 Next we present the bosonization of type $I$ and type $II$ 
vertex operators. For this purpose, 
let us consider the evaluation module. Set
\[ V_u=V\otimes_{\cA}\cA[u,u^{-1}],\quad
   V=\cA w_{+}\oplus \cA w_{-}. \]
We define the $ \cD Y_{\hbar}(\gtgl_2)$-module structure on $V_u$
as follows:
\begin{align*}
 H^{\sigma}(v).w_{\pm}=\frac{v-u\pm \hbar}{v-u}w_{\pm}\quad 
 \forall \sigma=\pm,&  \\
 K^{\pm}(v)\bigg\vert_{V_u}=\frac{v-u-\hbar}{v-u} \id_{V_u},
 \hspace{0.55 in} & \\
 E(v).w_{+}=0,\hspace{1.4 in} & F(v).w_{+}=\delta(v-u).w_{-}, \\
 E(v).w_{-}=\delta(v-u)w_{+},\hspace{0.8 in}& F(v).w_{-}=0. 
\end{align*}

\begin{defn}Vertex operators are intertwiners of the following form:
\begin{enumerate}
\renewcommand{\labelenumi}{(\roman{enumi})}
\item (type $I$) \hspace{1 in} 
      $\displaystyle{ \Phi^{(1-i,i)}(u):\cF_{i,s}\longrightarrow
                      \cF_{1-i,s-1}\otimes V_{u}, }$ 
\item (type $II$) \hspace{0.95 in}
      $\displaystyle{ \Psi^{(1-i,i)}(u):\cF_{i,s}\longrightarrow
                      V_{u}\otimes \cF_{1-i,s-1}. }$
\end{enumerate}
\end{defn}
Set 
\[ \Phi^{(1-i,i)}(u)=\sum_{\vep=\pm} \Phi^{(1-i.i)}_\vep(u)\otimes w_{\vep},
\qquad
   \Psi^{(1-i,i)}(u)=\sum_{\vep=\pm} w_{\vep}\otimes \Psi^{(1-i.i)}_\vep(u),
\]
We normalize them as 
\begin{enumerate}
\renewcommand{\labelenumi}{(\roman{enumi})}
\item $\displaystyle{ 
       \langle \Lambda_1,s-1|\Phi^{(1,0)}_{-}(u)|\Lambda_0,s \rangle=1,
       \qquad 
       \langle \Lambda_0,s-1|\Phi^{(0,1)}_{+}(u)|\Lambda_1,s \rangle=1,} $
\item $\displaystyle{
       \langle \Lambda_1,s-1|\Psi^{(1,0)}_{-}(u)|\Lambda_0,s \rangle=1,
       \qquad
       \langle \Lambda_0,s-1|\Psi^{(0,1)}_{+}(u)|\Lambda_1,s \rangle=1,} $
\end{enumerate}
where we set 
$|\Lambda_i,s\rangle=1\otimes
e^{\overline{\Lambda}_i+\frac{s}{2}(\vep_1+\vep_2)}$.
We mean by $\langle \Lambda_{1-i},s-1|
     \Phi^{(1-i,i)}_{\vep}(u)|\Lambda_{i},s \rangle$ the 
coefficient of $|\Lambda_{1-i},s-1 \rangle$ of the element
$\Phi^{(1-i,i)}_{\vep}(u)|\Lambda_{i},s \rangle$, and similarly for
$\Psi^{(1-i,i)}_{\vep}(u)$. With the above normalization our vertex 
operators uniquely exist. By using Lemma \ref{coproduct}, 
we obtain the bosonization formula of these vertex operators as follows.
\begin{thm}[Bosonization of vertex operators] \label{VOp}
\begin{align*}               
\Phi^{(1-i,i)}_{-}(u)&= \exp \left[ \sum_{k>0}
            \frac{a_{2,-k}}{k}(u+\frac{5}{4}\hbar )^k \right]
                        \exp \left[ \sum_{k>0}
            \frac{a_{1,k}}{k} (u+\frac{1}{4}\hbar )^{-k} \right] \\
                     &  \times e^{-\vep_2}
            \left[ -(u+\frac{1}{4}\hbar )\right]^{
            \frac{\partial_{\alpha}+i}{2}}, \\
\Phi^{(1-i,i)}_{+}(u)&= [\Phi^{(1-i,i)}_{-}(u), f_0], \\
\Psi^{(1-i,i)}_{+}(u)&= \exp \left[ \sum_{k>0}
            \frac{a_{1,-k}}{k}(u-\frac{1}{4}\hbar )^k \right]
                        \exp \left[ \sum_{k>0}
            \frac{a_{2,k}}{k} (u+\frac{3}{4}\hbar )^{-k} \right] \\
                     &  \times e^{-\vep_1}
            \left[ -(u+\frac{3}{4}\hbar )\right]^{
            \frac{-\partial_{\alpha}+i}{2}}(-1)^{1-i}, \\
\Psi^{(1-i,i)}_{-}(u)&= [\Psi^{(1-i,i)}_{+}(u), e_0].
\end{align*}
\end{thm}
Let us set
\begin{align*}
  \overset{1}{\Phi}^{(i,1-i)}(u_1)\overset{2}{\Phi}^{(1-i,i)}(u_2)
&=\sum_{\vep_1,\vep_2}
  \Phi^{(i,1-i)}_{\vep_1}(u_1)\Phi^{(1-i,i)}_{\vep_2}(u_2)
  w_{\vep_1}\otimes w_{\vep_2}, \\
  \overset{2}{\Phi}^{(i,1-i)}(u_2)\overset{1}{\Phi}^{(1-i,i)}(u_1)
&=\sum_{\vep_1,\vep_2}
  \Phi^{(i,1-i)}_{\vep_2}(u_2)\Phi^{(1-i,i)}_{\vep_1}(u_1)
  w_{\vep_1}\otimes w_{\vep_2}. 
\end{align*}
The commutation relations among vertex operators can be obtained from the
above theorem.
\begin{thm}[Commutation relation] \label{comm}
\begin{align*}
R(u_1-u_2)\overset{1}{\Phi}^{(i,1-i)}(u_1)\overset{2}{\Phi}^{(1-i,i)}(u_2)
&=\left( \frac{u_1+\frac{1}{4}\hbar}{u_2+\frac{1}{4}\hbar} \right)^{1-i}
  \overset{2}{\Phi}^{(i,1-i)}(u_2)\overset{1}{\Phi}^{(1-i,i)}(u_1), \\
R(u_1-u_2)\overset{1}{\Psi}^{(i,1-i)}(u_1)\overset{2}{\Psi}^{(1-i,i)}(u_2)
&=\left( \frac{u_1+\frac{3}{4}\hbar}{u_2+\frac{3}{4}\hbar} \right)^{i-1}
  \overset{2}{\Psi}^{(i,1-i)}(u_2)\overset{1}{\Psi}^{(1-i,i)}(u_1), \\
          \overset{1}{\Psi}^{(i,1-i)}(u_1)\overset{2}{\Phi}^{(1-i,i)}(u_2)
&=\left( \frac{u_1+\frac{3}{4}\hbar}{u_2+\frac{1}{4}\hbar} \right)^{1-i}
  \overset{2}{\Phi}^{(i,1-i)}(u_2)\overset{1}{\Psi}^{(1-i,i)}(u_1).
\end{align*}
\end{thm}
 We remark that all of the results obtained here are quite
natural deformation of that in $U_q(\hgtsl_2)$ case \cite{jimbo_miwa}
except for the grading operator.
\section{Discussion}\hspace{1 in} \\ 
In this paper, we have constructed a central extension of Yangian 
double $\cD Y_{\hbar}(\gtg)$ for $\gtg=\gtgl_2,\gtsl_2$.
We also introduced the Drinfel'd type generators which
are defined in \cite{dr_2} for the Yangian $Y_{\hbar}(\gtsl_2)$. 
We obtained the bosonization of level 1 
module and vertex operators of type $I$ and type $II$. Using explicit
expressions we also calculated their commutation relations.

One of the important applications of 
$\cD Y_{\hbar}(\gtsl_2)$ at level $0$ 
is F. Smirnov's theory \cite{smirnov} which states that 
the form factors in the Sine-Gordon theory satisfy 
deformed Knizhnik Zamolodchikov (d-KZ) equations. 
We expect that the correlation functions of the vertex 
operators satisfy the d-KZ equation at arbitrary level. 
For the lattice theory, like spin $\frac{1}{2}$ 
$XXX$-model, the space of state of the model might be constructed 
by level $1$ $ \cD Y_{\hbar}(\gtsl_2)$-module with its dual,
just as in the case of spin $\frac{1}{2}$ $XXZ$-model \cite{jimbo_miwa}. 

  From the mathematical point of view, it is interesting to
define the central extension of Yangian double 
$\cD Y_{\hbar}(\gtg)$ for an arbitrary simple finite dimensional Lie 
algebra $\gtg$. We would like to know the Chevalley generators of 
$\cD Y_{\hbar}(\gtg)$. To establish the representation theory of 
$\cD Y_{\hbar}(\gtg)$, especially highest weight modules, we need the 
triangular decomposition and the grading operator $d$. 
we should describe the universal $\cR$ of $\cD Y_{\hbar}(\gtg)$. 
We also have to establish an analogue of \cite{frenk_resh}'s theory. 
Another question which 
seems important is to construct free field realization of our algebra and
deformed Wakimoto-type module. 
 We hope to discuss these problems in our forthcoming papers.

\begin{ack}
The authors would like to thank J. Ding,
L. D. Faddeev, M. Jimbo, M. Kashiwara, T. Kohno, A. Leclair, 
A. Matsuo, A. Molev, T. Miwa, N. Yu. Reshetikhin, E. K. Sklyanin, 
and F. Smirnov for their interest and valuable suggestions. 
\end{ack}

\end{document}